\documentclass{amsart}

\usepackage{amsmath}
\usepackage{amstext}
\usepackage{amsfonts}
\usepackage{amscd}
\usepackage{amssymb,enumerate}
\usepackage{amsthm}
\usepackage{latexsym}
\usepackage[all]{xy}
\usepackage{setspace}

\newtheorem{lem}{Lemma}[section]
\newtheorem{cor}[lem]{Corollary}
\newtheorem{prop}[lem]{Proposition}
\newtheorem{thm}[lem]{Theorem}

\newtheorem{intthm}{Theorem}

\newtheorem{Defn}[lem]{Definition}
\newtheorem{Ex}[lem]{Example}
\newtheorem{Quest}[lem]{Question}
\newtheorem{Property}[lem]{Property}
\newtheorem{Properties}[lem]{Properties}
\newtheorem{Subprops}{}[lem]
\newtheorem{Para}[lem]{}
\newtheorem{Obs}[lem]{Observation}

\newtheorem{Remark}[lem]{Remark}

\newenvironment{para}{\begin{Para}\rm}{\end{Para}}

\theoremstyle{definition}

\newcommand{\ideal}[1]{\mathfrak{#1}}
\newcommand{\m}{\ideal{m}}

\newcommand{\p}{\ideal{p}}

\newcommand{\pcdim}{\operatorname{\PP_{\it C}\text{-}\pd}}

\newcommand{\icdim}{\operatorname{\I_{\it C}\text{-}\id}}

\newcommand{\pd}{\operatorname{pd}}

\newcommand{\id}{\operatorname{id}}

\newcommand{\im}{\operatorname{Im}}

\newcommand{\HH}{\operatorname{H}}
\newcommand{\Hom}{\operatorname{Hom}}

\newcommand{\coker}{\operatorname{Coker}}
\newcommand{\Ker}{\operatorname{Ker}}

\newcommand{\A}{\mathcal{A}}

\newcommand{\B}{\mathcal{B}}

\newcommand{\I}{\mathcal{I}}

\newcommand{\PP}{\mathcal{P}}

\newcommand{\ic}{\I_C}

\newcommand{\Rp}{R_{\p}}

\newcommand{\tot}{\operatorname{Tot}}

\newcommand{\tor}{\operatorname{Tor}}
\newcommand{\ext}{\operatorname{Ext}}

\newcommand{\xra}{\xrightarrow}

\renewcommand{\geq}{\geqslant}
\renewcommand{\ge}{\geqslant}
\renewcommand{\leq}{\leqslant}
\renewcommand{\le}{\leqslant}

\begin{document}
\dedicatory{Dedicated to the memory of Anders J. Frankild}

\author{Ryo Takahashi}\address{Department of Mathematical Sciences,
Faculty of Science, Shinshu University, 3-1-1 Asahi, Matsumoto,
Nagano 390-8621, Japan} \email{takahasi@math.shinshu-u.ac.jp}

\author{Diana White} \address{Department of Mathematics,
University of Colorado Denver, Campus Box 170, P.O. Box 173364
Denver, CO 80217-3364, USA}
   \email{diana.white@ucdenver.edu}

\thanks{This research was conducted in part while R.T. visited the University
of Nebraska in August 2006, partly supported by NSF grant DMS 0201904.}

\title{Homological aspects of semidualizing modules}


\begin{abstract}
We investigate the notion of the $C$-projective dimension of a
module, where $C$ is a semidualizing module.  When $C=R$, this
recovers the standard projective dimension.  We show that three
natural definitions of finite $C$-projective dimension agree, and
investigate the relationship between relative cohomology modules and
absolute cohomology modules in this setting.  Finally, we prove
several results that demonstrate the deep connections between
modules of finite projective dimension and modules of finite
$C$-projective dimension. In parallel, we develop the dual theory
for injective dimension and $C$-injective dimension.
\end{abstract}
\maketitle

\section*{Introduction}\label{sec:intro}

Grothendieck~\cite{hartshorne:lc} introduced dualizing modules as
tools for investigating cohomology theories in algebraic geometry.
In this paper, we investigate \emph{semidualizing} modules and
associated \emph{relative} cohomology functors.
Foxby~\cite{foxby:gmarm}, Vasconcelos~\cite{vasconcelos:dtmc} and
Golod~\cite{golod:gdgpi} independently initiated the study of
semidualizing modules. Over a noetherian ring $R$, a finitely
generated $R$-module $C$ is \emph{semidualizing} if the natural
homothety map $R \to \Hom_R(C,C)$ is an isomorphism and
$\ext^{\geqslant 1}_R(C,C)=0$; see~\ref{sdm} for a more general
definition. Examples include a dualizing module, when one exists,
and all finitely generated rank 1 projective modules.

Throughout this introduction, let $R$ be a commutative ring and $C$
a semidualizing $R$-module.  The class of \emph{$C$-projectives},
denoted $\PP_C$, consists of those $R$-modules of the form
$C\otimes_R P$ for some projective $R$-module $P$. These form the
building blocks of the so-called $G_C$-projectives, which are
studied in depth in~\cite{white:gcproj}.  Every $R$-module $M$
admits an \emph{augmented proper $\PP_C$-projective resolution}.
That is, there exists a complex
$$X^+=  \cdots \to C \otimes_R P_n\to \cdots \to C\otimes_R P_0 \to M \to 0$$
such that the complex $$\Hom_R(C,X^+)=  \cdots \to P_n \to \cdots
\to P_0 \to \Hom_R(C,M) \to 0$$ is exact. Despite the fact that
these augmented proper $\PP_C$-resolutions may not be exact, they
still have particularly good lifting properties.  In particular,
they give rise to well-defined cohomology modules
$\ext_{\PP_C}^i(M,N)$ for all $R$-modules $M$ and $N$;
see~\ref{proper}. 
 In the case $C=R$, these notions recover the
projectives, projective resolutions and the ``absolute'' cohomology
$\ext^i_R(M,N)$, respectively.

 Because augmented proper $\PP_C$-projective resolutions need not be exact, it is not
   immediately clear how best to define the $\PP_C$-projective dimension of a
   module.  For instance, should one consider arbitrary proper
$\PP_C$-projective resolutions
   or only exact ones?  Or should it be defined in terms of the vanishing
of the
   functors $\ext^n_{\PP_C}(M,-)$?  The next result, proved in Corollary~\ref{altdefn},
    shows that each of these
   approaches gives rise to the same invariant.

\begin{intthm}
Let $M$ be an $R$ module.  The following quantities are equal.
\begin{enumerate}[\rm(1)]
\item
 $\inf \left\{ \sup\{n\mid X_n\neq 0\}\left| \text{
\begin{tabular}{@{}c}
$X^+$ is an augmented \\
proper $\PP_C$-projective resolution of $M$
\end{tabular}
}\!\!\!\right. \right\}$

\

\item
 $\inf \left\{ \sup\{n\mid X_n\neq 0\}\left| \text{
\begin{tabular}{@{}c}
$X^+$ is an \emph{exact} augmented \\
proper $\PP_C$-projective resolution of $M$
\end{tabular}
}\!\!\!\right. \right\}$

\item $\sup\{n\mid \ext^n_{\PP_C}(M,-)\ne 0\}$
\end{enumerate}
\end{intthm}

The proof of this result uses so-called Auslander and Bass class
techniques; see~\ref{bass}.

Our investigation demonstrates a strong connection between the
modules of finite $\PP_C$-projective dimension and modules of finite
projective dimension, which is the focus of Section 2. For example,
the following is part of Theorem~\ref{either}.

\begin{intthm}If $M$ is an $R$-module, then
$\pcdim_R(M)=\pd_R(\Hom_R(C,M))$.
\end{intthm}

If $R$ is Cohen-Macaulay and local with dualizing module $D$, then a
result of Sharp shows that the modules of finite $\PP_D$-projective
dimension are precisely the modules of finite injective dimension.
Thus, this theorem recovers part of the Foxby equivalence
from~\cite{avramov:rhafgd}, namely that $M$ has finite injective
dimension if and only if $\Hom_R(D,M)$ has finite projective
dimension.

Section 4 explores the connection between the cohomology functors
$\ext_{\PP_C}$ and $\ext_R$. These connections are used
in~\cite{white:crct} to distinguish between several different
relative cohomology theories.  For example, the following is part of
Theorem~\ref{relabs} and Corollary~\ref{bcabs}. See~\ref{bass} for
the definition of the Bass class.

\begin{intthm}
Let $M$ and $N$ be $R$-modules.  There is an isomorphism
$\ext^i_{\PP_C}(M,N)\cong\ext^i_R(\Hom_R(C,M),\Hom_R(C,N))$ for all
$i$. If $M$ and $N$ are in the Bass class with respect to $C$, then
$\ext^i_{\PP_C}(M,N)\cong\ext^i_R(M,N)$ for all $i$.
\end{intthm}

This result and several others from this work have already proved
valuable for other investigations.  For instance, they are used by
the second author and her collaborators in~\cite{white:sgc} to
analyze the structure of certain categories naturally associated to
semidualizing modules, and in~\cite{white:crct,white:gcac} to study
balance properties of relative cohomology theories and to
distinguish between them.  We expect this line of inquiry to
continue to shed new light on the relationship between classical and
relative homological algebra.

Finally, in Section 5 of the paper we use results from the previous
sections to demonstrate the depth of the connection between modules
of finite $\PP_C$-projective dimension and modules of finite
projective dimension.

\section{Preliminaries}\label{sec:back}
Throughout this paper, let $R$ be a commutative ring.

  \begin{para}\label{complex}
  \label{cxs}
  An \emph{$R$-complex} is a sequence of $R$-module homomorphisms
  \begin{displaymath}
    X = \cdots\xra{\partial^X_{n+1}}X_n\xra{\partial^X_n}
    X_{n-1}\xra{\partial^X_{n-1}}\cdots
  \end{displaymath}
  such that $\partial^X_{n-1}\partial^X_{n}=0$ for each integer $n$; the
$n$th \emph{homology module} of $X$ is
$\HH_n(X)=\Ker(\partial^X_{n})/\im(\partial^X_{n+1})$. A morphism of
complexes $\alpha\colon X\to Y$ induces homomorphisms
$\HH_n(\alpha)\colon\HH_n(X)\to\HH_n(Y)$, and $\alpha$ is a
\emph{quasiisomorphism} when each $\HH_n(\alpha)$ is bijective. The
complex $X$ is \emph{bounded} if $X_n=0$ for $|n|\gg 0$.  It is
\emph{degreewise finite} if each $X_i$ is finitely generated.
\end{para}

\begin{para}\label{sdm}
An $R$-module $C$ is \emph{semidualizing} if
\begin{enumerate}[\quad(a)]
\item $C$ admits a degreewise finite projective resolution,
\item The natural homothety map $R \to \Hom_R(C,C)$ is an isomorphism, and
\item $\ext^{\geqslant 1}_R(C,C)=0$.
\end{enumerate}
A free $R$-module of rank 1 is semidualizing, and indeed this is the
only semidualizing module over a Gorenstein local ring. If $R$ is
noetherian, local and admits a dualizing module $D$, then $D$ is
semidualizing.
\end{para}

The next two classes of modules have been studied in numerous
papers, see e.g.~\cite{holm:cpaac} and~\cite{white:abc}.

\begin{para}\label{cproj}
The classes of \emph{$C$-projective} and \emph{$C$-injective}
modules are defined as
\begin{align*}
\PP_C & =\{\,C\otimes_R P\mid P\text{ is projective}\,\},\\
\ic & =\{\,\Hom_R(C,I) \mid I \text{ is injective}\,\}.
\end{align*}
When $C=R$, we omit the subscript and recover the classes of
projective and injective $R$-modules.
\end{para}

The next four paragraphs provide the necessary background on
relative homological algebra.  The reader is encouraged to
consult~\cite{enochs:rha} for details.

\begin{para}
The class $\PP_C$ is \emph{precovering} by~\cite[(5.10)]{white:abc}.
That is, given an $R$-module $M$, there exists a projective module
$P$ and a homomorphism $\phi\colon C\otimes_R P \to M$ such that,
for every projective $Q$, the induced map $$\Hom_R(C\otimes_R Q,
C\otimes_R P) \xra{\Hom_R(C\otimes_R Q, \phi)}  \Hom_R(C\otimes_R
Q,M)$$ is surjective.
  Dually, the class $\I_C$ is  \emph{preenveloping}.
\end{para}

\begin{para}\label{resolutions}
Since the class $\PP_C$ is precovering, for any $R$-module $M$ one
can iteratively take precovers to construct an \emph{augmented
proper $\PP_C$-projective resolution} of $M$, that is, a complex
\begin{displaymath}
    X^+ = \cdots\xra{\partial^X_{2}}C\otimes_R P_1
    \xra{\partial^X_{1}}C\otimes_R P_0 \xra{\partial^X_0} M\to 0
  \end{displaymath}
such that $\Hom_R(C\otimes_R Q,X^+)$ is exact for all projective
$R$-modules $Q$. The truncated complex
\begin{displaymath}
    X = \cdots\xra{\partial^{X^+}_{2}}C\otimes_R P_1
    \xra{\partial^{X^+}_{1}}C\otimes_R P_0\to 0
  \end{displaymath}
is a \emph{proper $\PP_C$-projective resolution} of $M$. For $n\ge
0$, set
$$
\Omega_n^{X^+}=
\begin{cases}
M & \text{if }n=0,\\
\Ker(\partial^{X^+}_{n-1}) & \text{if }n\ge 1.
\end{cases}
$$
Note that $X^+$ need not be exact unless $C=R$.

Dually, let $N$ be an $R$-module with \emph{augmented proper
$\I_C$-injective resolution}
$$
Y^+ = 0 \to N \to \Hom_R(C,I^0) \xra{\partial_Y^0} \Hom_R(C,I^1)
\xra{\partial_Y^1} \cdots
$$For an integer $n\ge 0$, set$$
\Omega_{Y^+}^{n}=
\begin{cases}
N & \text{if }n=0,\\
\im(\partial_{Y^+}^{n-1}) & \text{if }n\ge 1.
\end{cases}
$$

 \label{proper}
Proper $\PP_C$-projective resolutions are unique up to homotopy
equivalence; see e.g.~\cite[(1.8)]{holm:ghd}. Accordingly, the
\textit{$n$th relative cohomology modules}
\begin{align*}
\ext^n_{\PP_C}(M,N)&=\HH^n\Hom_R(X,N)
\end{align*}
where $X$ is a proper $\PP_C$-projective resolution of $M$ are
well-defined for each integer $n$.
The cohomology modules $\ext^n_{\I_C}(M,N)$ are defined dually.
\end{para}


\begin{para}  \label{xdim}
The \emph{$\PP_C$-projective dimension} of $M$
  is
  \begin{displaymath}
    \pcdim(M)= \inf \left\{ \sup\{n\mid X_n\neq 0\} \left|
      \begin{array}{l}
        \text{$X$ is a proper $\PP_C$-projective resolution of $M$}
      \end{array}
     \right. \!\!\!
     \right\}
  \end{displaymath}
  The modules of $\PP_C$-projective dimension zero are the non-zero
  modules in $\PP_C$.  The \emph{$\I_C$-injective dimension}, denoted
  $\icdim(-)$ is defined dually.
\end{para}

\begin{para}\label{dimshift}(\emph{Dimension Shifting}) Let $N$ be an
$R$-module.  For any augmented proper $\PP_C$-projective resolution
$X^+$ (as above) of $M$, there are isomorphisms
$$\ext^i_{\PP_C}(M,N)\cong \ext^{i-1}_{\PP_C}(\Omega_1^{X^+},N)
\cong \ext^{i-2}_{\PP_C}(\Omega_2^{X^+},N) \cong \cdots \cong
\ext^{i-n}_{\PP_C}(\Omega_n^{X^+},N)$$ for integers $1\leq n < i$.
\end{para}

\begin{para}\label{bass}The \textit{Bass class with respect to C},
 denoted $\B_C$ or $\B_C(R)$, consists of all $R$-modules $M$ satisfying
\begin{enumerate}[\quad\rm(a)]
\item\label{bass1} $\ext^{\geqslant 1}_R(C,M)=0=\tor_{\geqslant 1}^R(C,\Hom_R(C,M))=0$, and
\item\label{bass3} The natural evaluation map
 $\nu_{CCM}\colon C\otimes_R\Hom_R(C,M)\to M$ is an isomorphism.
We will write $\nu_M=\nu_{CCM}$ if there is no confusion.
\end{enumerate}
Dually, the \textit{Auslander class with respect to C},
 denoted $\A_C$ or $\A_C(R)$, consists of all $R$-modules $M$ satisfying
\begin{enumerate}[\quad\rm(a)]
\item\label{aus1} $\tor_{\geqslant 1}^R(C,M)=0=\ext_R^{\geq 1}(C,C\otimes_RM)$, and
\item\label{aus3} The natural map $\mu_{CCM}\colon
M\to\Hom_R(C,C\otimes_RM)$ is an isomorphism. We will write
$\mu_M=\mu_{CCM}$ if there is no confusion.
\end{enumerate}
\end{para}

We now state some basic results about the classes $\A_C$ and $\B_C$.
These facts are well-known when $R$ is noetherian. In this
generality, the first follows from~\cite[(6.5)]{white:abc}, the
second follows from the first, and the third assertion is routine to
check.

\begin{para}\label{projac}The following hold.
\begin{enumerate}[\quad\rm (a)]
\item\label{twothree} If any two $R$-modules in a short exact
sequence are in $\A_C$, respectively $\B_C$, then so is the third.
\item\label{pdimac} \label{bcinj} The class $\A_C$ contains all modules of
finite flat dimension.  The class $\B_C$ contains all modules of
finite injective dimension.
\item\label{acbccorr}If $M$ is in $\A_C$, then $C\otimes_R M$ is
in $\B_C$.  If $M$ is in $\B_C$, then $\Hom_R(C,M)$ is in $\A_C$.
\end{enumerate}
\end{para}

\section{Relative dimensions and Auslander and Bass classes}

This section has two interwoven themes.  First, we explore the
interplay between $\PP_C$-projective dimension and projective
dimension. Second, we investigate the exactness of augmented proper
$\PP_C$-projective resolutions for modules of finite
$\PP_C$-projective dimension.

We begin with the following lemma, which follows from the
definitions of semidualizing modules and augmented proper
resolutions, using~\ref{projac}\eqref{bcinj}.

\begin{lem}\label{cdimpdim}\label{pcpdim}
Let $C$ be a semidualizing $R$-module and $M$ an $R$-module.
\begin{enumerate}[\quad\rm(a)]
\item If $X^+$ is an augmented proper $\PP_C$-projective resolution of $M$,
then $\Hom_R(C,X^+)$ is an augmented projective resolution of
$\Hom_R(C,M)$.
\item If $Y^+$ is an augmented proper $\I_C$-injective resolution
of $M$, then $C\otimes_R Y^+$ is an augmented injective resolution
of $C\otimes_R M$.\qed
\end{enumerate}
\end{lem}

We now investigate exactness of augmented proper $\PP_C$-projective
resolutions.

\begin{prop}\label{somevery}
 Let $C$ be a semidualizing $R$-module, $M$ an $R$-module and $n$ a
 nonnegative integer.
\begin{enumerate}[\quad\rm(a)]
\item The following are equivalent.
\begin{enumerate}[\quad\rm(i)]
\item There exists an augmented proper $\PP_C$-projective resolution of $M$
 which is exact in degree less than $n$;
\item Every augmented proper $\PP_C$-projective resolution of $M$ is exact
in degree
less than $n$;
\item The natural homomorphism $\nu_M:C\otimes_R\Hom_R(C,M)\to M$ is an
isomorphism and $\tor^R_i(C,\Hom_R(C,M))=0$ for $0<i<n$.
\end{enumerate}
\item The following are equivalent.
\begin{enumerate}[\quad\rm(i)]
\item There exists an augmented proper $\I_C$-injective resolution of $M$
which is
exact in degree less than $n$;
\item All augmented proper $\I_C$-injective resolutions of $M$ are exact in
degree less than $n$;
\item The natural homomorphism $\mu_M:C\to\Hom_R(C,C\otimes_RM)$ is an
isomorphism and $\ext_R^i(C,C\otimes_RM)=0$ for $0<i<n$.
\end{enumerate}
\end{enumerate}
\end{prop}

\begin{proof}We prove only part(a).
Let $X^+$ be an augmented proper $\PP_C$-projective resolution of
$M$. By Lemma~\ref{cdimpdim}, the complex $\Hom_R(C,X)$ is a
projective resolution of $\Hom_R(C,M)$. Hence, there is an
isomorphism $\tor_i^R(C,\Hom_R(C,M))\cong H_i(C\otimes_R\Hom_R(C,
X))$
 for all $i\ge 0$.
Since each $X_i$ is in $\B_C$, the natural chain map
$C\otimes_R\Hom_R(C, X)\to X$ is an isomorphism.  The result follows
\end{proof}

From the above proposition, we obtain the following criterion for a
given module to possess exact augmented proper resolutions.

\begin{cor}\label{exactpcres}Let $C$ be a semidualizing
$R$-module and $M$ an $R$-module.
\begin{enumerate}[\quad\rm(a)]
\item The following are equivalent.
\begin{enumerate}[\quad\rm(i)]
\item $M$ admits an exact augmented proper $\PP_C$-projective resolution;
\item All augmented proper $\PP_C$-projective resolutions of $M$ are exact;
\item The natural homomorphism $\nu_{M}\colon C\otimes_R\Hom_R(C,M)\to M$ is an
isomorphism and $\tor^R_{\geq 1}(C,\Hom_R(C,M))=0$.
\end{enumerate}
\item The following are equivalent.
\begin{enumerate}[\quad\rm(i)]
\item $M$ admits an exact augmented proper $\I_C$-injective resolution;
\item All augmented proper $\I_C$-injective resolutions of $M$ are exact;
\item The natural homomorphism $\mu_{M}\colon M\to\Hom_R(C,C\otimes_RM)$ is an
isomorphism and $\ext_R^{\geq 1}(C,C\otimes_RM)=0$.\qed
\end{enumerate}
\end{enumerate}
\end{cor}

From the definitions of the Auslander and Bass classes, we have the
following.

\begin{cor}\label{bcpcres}Let $C$ be a semidualizing $R$-module and $M$ an
$R$-module.
\begin{enumerate}[\quad\rm (a)]
\item Assume $M$ is in $\B_C$.  Then every augmented proper $\PP_C$-projective
resolution of $M$ is then exact. In particular, every
$\PP_C$-precover of $M$ is surjective.
\item Assume $M$ is in $\A_C$. Then every augmented proper $\I_C$-resolution of $M$
is exact. In particular, every $\I_C$-preenvelope of $M$ is
injective.\qed
\end{enumerate}
\end{cor}

The next few technical results build toward the following fact,
which follows from Corollaries~\ref{exactpcres} and~\ref{altdefn}:
if $M$ has finite $\PP_C$-projective dimension, then \emph{every}
augmented proper $\PP_C$-resolution of $M$ is exact.

\begin{lem}\label{identity}
Let $C$ be a semidualizing module and $M$ an $R$-module.
\begin{enumerate}[\quad\rm(a)]
\item
The composition $\Hom_R(C,\nu_M)\circ\mu_{\Hom_R(C,M)}$ is the
identity map on $\Hom_R(C,M)$. Hence, $\Hom_R(C,\nu_M)$ is a split
epimorphism and $\mu_{\Hom_R(C,M)}$ is a split monomorphism.
\item
Assume that $\nu_M$ is injective. The composition
$\mu_{\Hom_R(C,M)}\circ\Hom_R(C,\nu_M)$ is the identity map on
$\Hom_R(C, C\otimes_R\Hom_R(C,M))$. Hence, $\Hom_R(C,\nu_M)$ is an
isomorphism and $\mu_{\Hom_R(C,M)}$ is the inverse isomorphism.
\item
The composition $\nu_{C\otimes_RM}\circ(C\otimes_R\mu_M)$ is the
identity map on $C\otimes_R M$. Hence, $C\otimes_R\mu_M$ is a split
monomorphism and $\nu_{C\otimes_RM}$ is a split epimorphism.
\item
Assume that $\mu_M$ is surjective. The composition
$(C\otimes_R\mu_M)\circ\nu_{C\otimes_RM}$ is the identity map on
$C\otimes_R\Hom_R(C,C\otimes_R M)$. Hence, $C\otimes_R\mu_M$ is an
isomorphism and $\nu_{C\otimes_RM}$ is the inverse isomorphism.
\end{enumerate}
\end{lem}

\begin{proof}Part (a) is straightforward to check.  For part (b),
set $\rho =\mu_{\Hom_R(C,M)}\circ\Hom_R(C,\nu_M)$.  Note that if
$\xi\in\Hom_R(C,C\otimes_R\Hom_R(C,M))$, then $\rho(\xi)$ sends
$z\in C$ to $z\otimes(\nu_M\circ\xi)$. Thus,
$\nu_M(\xi(z)-z\otimes(\nu_M\cdot\xi))=\nu_M(\xi(z))-\nu_M(\xi(z))=0$.
Since $\nu_M$ is injective, there is an equality
$\xi(z)=z\otimes(\nu_M\cdot\xi)$, thereby showing that $\rho$ is the
identity. Parts (c) and (d) are proved similarly.
\end{proof}

\begin{lem}\label{general}Let $C$ be a semidualizing $R$-module and
$M$ an $R$-module.
\begin{enumerate}[\quad\rm(a)]
\item Assume $\ext_R^1(C,C\otimes_R\Hom_R(C,M))=0$. If $\nu_M$ is injective,
then it is an isomorphism.
\item Assume that $\tor_1^R(C,\Hom_R(C,C\otimes_RM))=0$. If $\mu_M$ is
surjective, then it is an isomorphism.
\end{enumerate}
\end{lem}

\begin{proof}We prove only part (a), as part (b) is dual. Set $L=\coker{\nu_M}$.
Applying $\Hom_R(C,-)$ to the exact sequence $$0 \to
C\otimes_R\Hom_R(C,M) \overset{\nu_M}{\longrightarrow} M \to L \to
0$$ induces an exact sequence
$$0 \to \Hom_R(C,C\otimes\Hom_R(C,M)) \xra{\Hom(C,\nu_M)}\\
\Hom_R(C,M) \to \Hom_R(C,L) \to 0$$ where right-exactness follows
from the equality $\ext_R^1(C,C\otimes_R\Hom_R(C,M))=0$. By
Lemma~\ref{identity}, $\Hom_R(C,\nu_M)$ is an isomorphism. Hence,
the above exact sequence implies that $\Hom_R(C,L)=0$, and it
follows from~\cite[(3.6)]{white:abc} that $L=0$.
\end{proof}

\begin{cor}\label{injbij}Let $C$ be a semidualizing $R$-module and $M$ an
 $R$-module.
\begin{enumerate}[\quad\rm (a)]
\item Suppose $\Hom_R(C,M)$ is in $\A_C$. If $\nu_M$ is injective,
then it is an isomorphism.
\item Suppose $C\otimes_RM$ is in $\B_C$. If $\mu_M$ is surjective,
then it is an isomorphism.\qed
\end{enumerate}
\end{cor}

We now prove our first theorem, which leads to some of the main
results of this section.

\begin{thm}\label{remove}
Let $C$ be a semidualizing $R$-module and $M$ an $R$-module.
Then the following hold.
  \begin{enumerate}[\quad\rm(a)]
 \item $M\in \B_C$ if and only if $\Hom_R(C,M)\in\A_C$.
 \item $M\in \A_C$ if and only if $C\otimes_R M\in\B_C$.
  \end{enumerate}
\end{thm}

\begin{proof}We prove only part (a), as part (b) is proved
similarly. According to \ref{projac}\eqref{acbccorr}, it is enough
to assume $\Hom_R(C,M)\in\A_C$ and show $M\in \B_C$. The definition
of $\A_C$ implies that $\tor^R_{\geq 1}(C,\Hom_R(C,M))=0$ and
$\ext^{\geq 1}_R(C,C\otimes_R\Hom_R(C,M))=0$.  We will show that
that the evaluation map $\nu_M$ is an isomorphism.  Using the above
vanishings, it will then follow that $\ext^{\geq 1}_R(C,M)=0$ and
$M$ is in $\B_C$.

By Lemma~\ref{identity}(a), the composition
$\Hom_R(C,\nu_M)\circ\mu_{\Hom_R(C,M)}$ is the identity map on
$\Hom_R(C,M)$.  Since $\Hom_R(C,M)$ is in $\A_C$, the map
$\mu_{\Hom_R(C,M)}$ is an isomorphism, and so $\Hom_R(C,\nu_M)$ is
also an isomorphism. Setting $K=\ker(\nu_M)$, it follows that
$\Hom_R(C,K)=0$, and so by~\cite[(3.6)]{white:abc}, $K=0$.  Thus,
the map $\nu_M$ is injective, hence an isomorphism by
Corollary~\ref{injbij}(a).
\end{proof}

Recall that the class $\B_C$ contains all modules of finite
injective dimension and the class $\A_C$ contains all modules of
finite projective dimension; see~\ref{projac}.  By virtue of Theorem
\ref{remove}, we now obtain additional examples of modules in $\A_C$
and $\B_C$.

\begin{cor}\label{pcbc}Let $C$ be a semidualizing $R$-module and $M$ an
$R$-module.
\begin{enumerate}[\quad\rm (a)]
\item If $\pcdim_R(M)$ is finite, then $M$ is in $\B_C$.
\item If $\icdim_R(M)$ is finite, then $M$ is in $\A_C$.
\end{enumerate}
\end{cor}

\begin{proof}
We prove only part (a). Assume $\pcdim_R(M)$ is finite, and let
$X^+$ be an augmented proper $\PP_C$-resolution of $M$.  By
Lemma~\ref{cdimpdim}(a), the complex $\Hom_R(C,X^+)$ is a bounded
projective resolution of $\Hom_R(C,M)$.
 By~\ref{projac}\eqref{pdimac}, $\Hom_R(C,M)$ is in $\A_C$
 so Theorem~\ref{remove}(a) implies that $M$ is in $\B_C$.
\end{proof}

This yields the following key result.

\begin{cor}\label{altdefn}
Let $C$ be a semidualizing $R$-module and $M$ an $R$-module.
\begin{enumerate}[\quad\rm(a)]
\item The inequality $\pcdim_R(M)\leq n$ holds if and only if there
is an exact sequence
$$0 \to C\otimes_R P_n \to \cdots \to C\otimes_R P_0 \to M \to 0$$
with each $P_i$ a projective $R$-module.
\item  The inequality $\icdim_R(M)\leq n$ holds if and only if there
is an exact sequence
$$0 \to M \to \Hom_R(C,I^0) \to \cdots \to \Hom_R(C,I^n) \to 0$$
with each $I^i$ an injective $R$-module.
\end{enumerate}
\end{cor}

\begin{proof}
The ``only if'' direction of each statement follows immediately from
Corollaries \ref{bcpcres} and \ref{pcbc}.  For the other
implications,
 use dimension shifting~\ref{dimshift} to see that any exact sequence
 of the given
form is an augmented proper $\PP_C$-projective or $\I_C$-injective
resolution of $M$.
\end{proof}

We now investigate how $\PP_C$-projective and projective dimension
relate.

\begin{thm}\label{dimension}\label{either}Let $C$ be a semidualizing $R$-module.
The following equalities hold.
\begin{enumerate}[\quad\rm(a)]
\item $\pd_R(M)=\pcdim_R(C\otimes_R M)$
\item $\icdim_R(M)=\id_R(C\otimes_R M)$
\item $\pcdim_R(M)=\pd_R(\Hom_R(C,M))$
\item $\id_R(M) =\icdim_R(\Hom_R(C,M))$
\end{enumerate}
\end{thm}

\begin{proof}
We prove only part (a). Assume $\pd_R(M)=s<\infty$ and consider an
augmented projective resolution of M
$$X= \quad 0 \to P_s \to P_{s-1} \to \cdots \to P_0 \to M \to 0.$$
 By~\ref{projac}(b), one has $M\in\A_C$ and so
$\tor_{\geq 1}^R(C,M)=0$.  Thus, the complex
$$
C\otimes_R X= \quad 0 \to C\otimes_RP_s \to C\otimes_RP_{t-1} \to
\cdots \to C\otimes_RP_0 \to C\otimes_RM \to 0
$$
is exact and thus an augmented proper $\PP_C$-projective resolution
of $C\otimes_RM$.  Note that properness can be shown by
using~\ref{dimshift}, or as a special case
of~\cite[(4.4)]{white:gcproj}. This provides an inequality
$s\ge\pcdim_R(C\otimes_RM)$. Conversely, assume that
$\pcdim_R(C\otimes_RM)=t<\infty$. By Corollary \ref{altdefn}(a),
there is an augmented \emph{exact} proper $\PP_C$-resolution of
$C\otimes_RM$
$$
X^+=\quad 0 \to C\otimes_RP_t \to C\otimes_RP_{t-1} \to \cdots \to
C\otimes_RP_0 \to C\otimes_RM \to 0.
$$
Thus, the complex $\Hom_R(C,X^+)$ is exact. Corollary~\ref{pcbc}
implies that $C\otimes_RM$ is in $\B_C$ and Theorem~\ref{remove}
then implies that $M$ is in $\A_C$.  Thus, $\mu_M$ is an
isomorphism.  Since each $\mu_{P_i}$ is also an isomorphism,
$\Hom_R(C,X^+)$ is isomorphic to an exact sequence of the form $$0
\to P_t \to P_{t-1} \to \cdots \to P_0 \to M \to 0.$$ Thus,
$\pd_R(M)\le t=\pcdim_R(C\otimes_RM)$.
\end{proof}

Note that by assembling the information above, we get the following
extension of the Foxby equivalence~\cite{avramov:rhafgd}.

\begin{thm}\label{foxby}(Foxby equivalence) Let $C$ be a semidualizing
$R$-module, and let $n$ be a non-negative integer.  Set
$\widehat{\PP_C}(R)_{\leq n}$, $\widehat{\PP}(R)_{\leq n}$,
$\widehat{\I_C}(R)_{\leq n}$, and $\widehat{\I}(R)_{\leq n}$ to be
the classes of modules of $C$-projective, projective, $C$-injective
and injective dimension of at most $n$, respectively. Then there are
equivalences of categories
  \begin{displaymath}
    \xymatrix@C=20ex{\PP(R) \ar@{^(->}[d] \ar@<0.8ex>[r]_-{\sim}^-{C
        \otimes_R-} & \PP_C(R)\ar@{^(->}[d]
      \ar@<0.8ex>[l]^-{\Hom_R(C,-)} \\
      \widehat{\PP}(R)_{\leq n} \ar@{^(->}[d] \ar@<0.8ex>[r]_-{\sim}^-{C
        \otimes_R-} & \widehat{\PP}_C(R)_{\leq n} \ar@{^(->}[d]
      \ar@<0.8ex>[l]^-{\Hom_R(C,-)} \\
      \A_C(R) \ar@<0.8ex>[r]_-{\sim}^-{C \otimes_R-} & \B_C(R)
      \ar@<0.8ex>[l]^-{\Hom_R(C,-)} \\
      \widehat{\I}_C(R)_{\leq n} \ar@{^(->}[u] \ar@<0.8ex>[r]_-{\sim}^-{C \otimes_R-} &
      \widehat{\I}(R)_{\leq n}. \ar@<0.8ex>[l]^-{\Hom_R(C,-)} \ar@{^(->}[u] \\
       \I_C(R) \ar@{^(->}[u] \ar@<0.8ex>[r]_-{\sim}^-{C \otimes_R-} &
      \I(R). \ar@<0.8ex>[l]^-{\Hom_R(C,-)} \ar@{^(->}[u]
    }
  \end{displaymath}
\end{thm}


\section{Vanishing of relative cohomology and consequences}

In this section, we investigate how vanishing of the relative
cohomology functors $\ext^i_{\PP_C}(M,-)$ and $\ext^i_{\I_C}(-,N)$,
respectively, characterizes the finiteness of $\pcdim_R(M)$ and
$\icdim_R(N)$.  Our proofs use specific properties of semidualizing
modules and so do not carry over directly to other relative
settings.

\begin{thm}\label{charcproj}Let $C$ be a semidualizing
$R$-module and $M$ an $R$-module.
\begin{enumerate}[\quad\rm (a)]
\item The following are equivalent.
\begin{enumerate}[\quad\rm (i)]
\item $\ext^1_{\PP_C}(M,-)=0$
\item $\ext^{\geq 1}_{\PP_C}(M,-)=0$
\item $M$ is $C$-projective
\end{enumerate}
\item The following are equivalent.
\begin{enumerate}[\quad\rm (i)]
\item $\ext^1_{\I_C}(-,M)=0$
\item $\ext^{\geq 1}_{\I_C}(-,M)=0$
\item $M$ is $C$-injective\
\end{enumerate}
\end{enumerate}
\end{thm}

\begin{proof}
We prove part (a); part (b) is dual.

(iii) $\Rightarrow$ (ii): If $M$ is $C$-projective, then the complex
$$\cdots \to 0 \to M \overset{=}{\to} M \to 0$$
is an augmented proper $\PP_C$-resolution of $M$ and so
$\ext_{\PP_C}^{\geq 1}(M,-)=0$.

(ii) $\Rightarrow$ (i) is immediate.

(i) $\Rightarrow$ (iii): Let
$$
X= \cdots \overset{d_2}{\to} C\otimes_RP_1 \overset{d_1}{\to}
C\otimes_RP_0 \xra{d_0} M \to 0$$ be an augmented proper
$\PP_C$-resolution of $M$. Let $K_0$ be the kernel of $d_0$, and let
$\beta\colon K_0\to C\otimes_RP_0$ be the inclusion map. There is a
homomorphism $\alpha\colon C\otimes_RP_1\to K_0$ such that
$d_1=\beta\alpha$. Noting that $\beta\alpha d_2=d_1d_2=0$, the
injectivity of $\beta$ implies that $\alpha d_2=0$. Since
$\ext_{\PP_C}^1(M,K_0)=0$, the induced sequence
$$
\Hom_R(C\otimes_RP_0,K_0) \to \Hom_R(C\otimes_RP_1,K_0) \to
\Hom_R(C\otimes_RP_2,K_0)
$$
is exact.  Hence, there exists $\xi\in\Hom_R(C\otimes_RP_0,K_0)$
such that $\alpha=\xi d_1=\xi\beta\alpha$. There is an equality
$\Hom_R(C,\alpha)=\Hom_R(C,\xi)\circ\Hom_R(C,\beta)\circ
\Hom_R(C,\alpha)$ and so
$\Hom_R(C,\xi)\circ\Hom_R(C,\beta)=\id_{\Hom_R(C,K_0)}$, as
$\Hom_R(C,\alpha)$ is surjective. Therefore, the exact sequence
(which is exact by properness of $X$)
$$
0 \to \Hom_R(C,K_0) \xra{\Hom_R(C,\beta)} \Hom_R(C,C\otimes_RP_0)
\to \Hom_R(C,M) \to 0
$$
splits. Since $\Hom_R(C,C\otimes_RP_0)\cong P_0$ is $R$-projective,
so is $\Hom_R(C,M)$. Theorem~\ref{either} and
Corollary~\ref{altdefn}(a) now imply that $M$ is $C$-projective.
\end{proof}

Using dimension shifting, we have the following extension of the
previous result.

\begin{thm}\label{five}Let $C$ be a semidualizing $R$-module, let
$n$ be a non-negative integer, and let $M,N$ be $R$-modules.
\begin{enumerate}[\quad\rm (a)]
\item The following are equivalent.
\begin{enumerate}[\quad\rm(i)]
\item $\ext^{n+1}_{\PP_C}(M,-)=0$
\item $\ext^{\ge n+1}_{\PP_C}(M,-)=0$
\item $\pcdim_R(M)\le n$.
\end{enumerate}
\item The following are equivalent.
\begin{enumerate}[\quad\rm(i)]
\item $\ext^{n+1}_{\I_C}(-,N)=0$
\item $\ext^{\ge n+1}_{\I_C}(-,N)=0$
\item $\icdim_R(N)\le n$.\qed
\end{enumerate}
\end{enumerate}
\end{thm}

The preceding two results imply the following.

\begin{cor}\label{schanuel} Let $C$ be a semidualizing
$R$-module, $M$ an $R$-module and $n\geq 0$.
\begin{enumerate}[\quad\rm(a)]
\item The following are equivalent.
\begin{enumerate}[\quad\rm(i)]
\item There exists an augmented proper
$\PP_C$-projective resolution $X^+$ of $M$ such that
$\Omega_n^{X^+}$ is $C$-projective.
\item Every augmented proper
$\PP_C$-projective resolution $X^+$ of $M$ has the property that
$\Omega_n^{X^+}$ is $C$-projective.
\end{enumerate}
\item The following are equivalent.
\begin{enumerate}[\quad\rm(i)]
\item There exists an augmented proper
$\I_C$-injective resolution $Y^+$ of $M$ such that $\Omega_{Y^+}^n$
is $C$-injective.
\item  Every augmented proper
$\I_C$-injective resolution $Y^+$ of $M$ has the property that
$\Omega_{Y^+}^n$ is $C$-injective.\qed
\end{enumerate}
\end{enumerate}
\end{cor}

We conclude this section by showing that if any two of three modules
in a short exact sequence have finite $\PP_C$-projective dimension
then so does the third.  Note that the standard constructions (using
Horseshoe Lemmas, mapping cones, etc.) that show this result when
$C=R$ can be used in this setting. However, we offer a shorter proof
that uses the classical result.

\begin{prop}\label{cdim23}
Let $C$ be a semidualizing $R$-module.  Consider an exact sequence
of $R$-modules
$$
0\to M' \to M \to M'' \to 0.
$$
If any two of the modules have finite $\PP_C$-projective dimension,
respectively
 $\I_C$-injective dimension, then so does the third.
\end{prop}

\begin{proof}
Assume that two of the modules $M,M',M''$ have finite
$\PP_C$-projective dimension. By Corollary~\ref{pcbc}, these two
modules are in $\B_C$. By~\ref{projac}\eqref{twothree}, this forces
all of the modules $M,M',M''$ to be in $\B_C$. Thus, the complex
$$
\ \ 0 \to \Hom_R(C,M') \to \Hom_R(C,M) \to \Hom_R(C,M'') \to 0
$$
is exact. By Theorem~\ref{dimension}(c), two of the above
  Hom modules have finite projective dimension and hence so
 does the third.
Theorem~\ref{either}(c) implies that all of  $M,M',M''$
 have finite $P_C$-projective dimension.  The other assertion is
 dual.
\end{proof}

\section{Comparing relative and absolute cohomology}

In this section we investigate the interplay between absolute Ext
and the relative cohomologies $\ext_{\PP_C}$ and $\ext_{\I_C}$.
Under certain circumstances, they agree with their corresponding
absolute counterparts.

\begin{thm}\label{relabs}
Let $C$ be a semidualizing $R$-module, and let $M$ and $N$ be
$R$-modules.  There exist isomorphisms.
\begin{align*}
\ext^i_{\PP_C}(M,N)&\cong\ext^i_R(\Hom_R(C,M),\Hom_R(C,N))\\
\ext^i_{\I_C}(M,N)&\cong\ext^i_R(C\otimes_R M,C\otimes_R N)
\end{align*}
\end{thm}

\begin{proof} We prove only part (a).
Let $X^+$ be an augmented proper $\PP_C$-resolution of $M$. By
Lemma~\ref{cdimpdim}, the complex $\Hom_R(C,X^+)$ is an augmented
projective resolution of $\Hom_R(C,M)$. Thus, the equalities below
hold by definition
\begin{align*}
\ext_R^i(\Hom_R(C,M),\Hom_R(C,N)) & = \HH^i(\Hom_R(\Hom_R(C,X),\Hom_R(C,N))) \\
& \cong \HH^i(\Hom_R(C\otimes_R\Hom_R(C,X),N)) \\
& \cong \HH^i(\Hom_R(X,N)) \\
& = \ext_{\PP_C}^i(M,N),
\end{align*}
while the isomorphisms follow from adjunction and the containment
$\PP_C \subseteq\B_C$.
\end{proof}

With appropriate Auslander and Bass class assumptions, the
aforementioned relative cohomology modules agree precisely with the
absolute Ext.

\begin{cor}\label{bcabs}
Let $C$ be a semidualizing $R$-module, and let $M,N$ be $R$-modules.
\begin{enumerate}[\quad\rm(a)]
\item
If $M$ and $N$ are in $\B_C$, then $\ext^i_{\PP_C}(M,N)\cong\ext^i_R(M,N)$ for all $i$.
\item
If $M$ and $N$ are in $\A_C$, then $\ext^i_{\I_C}(M,N)\cong\ext^i_R(M,N)$ for all $i$.
\end{enumerate}
\end{cor}

\begin{proof}We prove only part(a).  Let $\tot X$ denote the total
complex of a double complex $X$.  Let $P^+$ be an augmented
projective resolution of $\Hom_R(C,M)$, and let $I^+$ be an
augmented injective resolution of $N$. Since $M$ and $N$ are in
$\B_C$, the complexes $C\otimes_RP^+$ and $\Hom_R(C,I^+)$ are exact.
There is an isomorphism
$$
\Hom_R(C\otimes_RP,I)\cong\Hom_R(P,\Hom_R(C,I))
$$
of double complexes. This provides the second isomorphism below
\begin{align*}
\ext_R^i(M,N) & \cong H^i(\tot\Hom_R(C\otimes_RP,I)) \\
& \cong H^i(\tot\Hom_R(P,\Hom_R(C,I))) \\
& \cong \ext_R^i(\Hom_R(C,M),\Hom_R(C,N)) \\
& \cong \ext_{\PP_C}^i(M,N),
\end{align*}
while the last isomorphism follows from Theorem~\ref{relabs}.
\end{proof}

\section{Further parallels between the classical and relative
theories}

  The results of the previous sections demonstrate that there is a
   tight connection between modules of finite $\PP_C$-projective
   dimension and modules of finite projective dimension.  In this
   section we indicate how the machinery developed above allows
   us to extend many classical results to this new setting.  We begin
   by showing that $\PP_C$-projective dimension has the ability to detect when
   a ring is regular.

\begin{prop}\label{abs}
Let $(R,\m,k)$ be a noetherian, local ring and $C$ a semidualizing
$R$-module. The following are equivalent.
\begin{enumerate}[\quad\rm(i)]
\item\label{abs1} $\pcdim_R(M)$ is finite for all $R$-modules $M$.
\item\label{abs2} $\pcdim_R(k)$ is finite.
\item\label{abs3} $R$ is regular.
\end{enumerate}
\end{prop}

\begin{proof}
$\eqref{abs1}\implies\eqref{abs2}$ is trivial.

$\eqref{abs2}\implies\eqref{abs3}$ Since $\pcdim_R(k)$ is finite,
Lemma~\ref{cdimpdim} implies $\pd_R(\Hom_R(C,k))$ is finite. Since
$\Hom_R(C,k)$ is a nonzero $k$-vector space, $\pd_R(k)$ is finite.
Thus, $R$ is regular.

$\eqref{abs3}\implies\eqref{abs1}$ Since $R$ is regular, the only
semidualizing $R$-module is $R$ itself.  Thus, $C=R$ so this follows
from the Auslander-Buchsbaum-Serre theorem.
\end{proof}

Our methods also apply to bounded complexes of $C$-projective
modules, as the next result shows.

\begin{prop}[New Intersection Theorem for complexes of $C$-projective modules]
Let $(R,\m)$ be a noetherian local ring and $C$ a semidualizing
$R$-module.  If there exists a non-exact complex $$X= 0 \to
C^{\alpha_s} \to C^{\alpha_{s-1}} \to \cdots \to C^{\alpha_1} \to
C^{\alpha_0} \to 0$$ with $\ell_R(\HH_i(X))$ finite for all $i$,
then $s \geq \dim(R)$.
\end{prop}

\begin{proof}
First, note that the complex $$\Hom_R(C,X)= 0\to R^{\alpha_s} \to
R^{\alpha_{s-1}} \to \cdots \to R^{\alpha_1} \to R^{\alpha_0} \to
0$$ is non-exact.  Indeed, if it were exact, then it would split,
forcing the complex $C\otimes_R\Hom_R(C,X)\cong X$ to be exact, a
contradiction.

Now fix a prime $\p\neq\m$.  Since the homology of the complex $X$
has support equal to $\m$, the complex $X_{\p}$ is exact.  This
forces the complex $\Hom_{\Rp}(C_{\p},X_{\p})$ to be exact, as
$\ext^{\geqslant 1}_R(C,C)=0$. This forces
$\HH_i(\Hom_R(C,X))_{\p}\cong\HH_i(\Hom_{\Rp}(C_{\p},X_{\p}))=0$ for
all $i$. Thus, $\ell_R(H_i(\Hom_R(C,X)))<\infty$ for all $i$.  The
New Intersection Theorem now implies that $s\geqslant \dim(R)$.
\end{proof}

Next, we extend Bass' result that a ring is noetherian if and only
if the class of injective $R$-modules is closed under direct sums.

\begin{prop}
Let $R$ be a commutative ring and $C$ a semidualizing $R$-module.
 The ring $R$ is noetherian if and only if the class $\I_C$ is closed
  under direct sums.
\end{prop}

\begin{proof}
Assume $R$ is noetherian.  Let $\{I_{\lambda}\}$ be a collection of
injective $R$-modules. Since $C$ is finitely presented, there is an
isomorphism
$$\coprod_{\lambda}\Hom_R(C,I_{\lambda})\cong\Hom_R(C,
\coprod_{\lambda}I_{\lambda})$$ and the desired result follows by
Bass' result.

Conversely, assume the class $\I_C$ is closed under direct sums.
 Let $\{I_{\lambda}\}$ be a collection of injective modules so that
 $\{\Hom_R(C, I_{\lambda})\}$ is a collection of $C$-injective
modules. By assumption, there is an isomorphism
$$\coprod \Hom_R(C,I_{\lambda})\cong\Hom_R(C,J)$$ for some injective $J$.
This provides the third isomorphism below. The first and fourth isomorphisms follow from the fact that any injective module is in $\B_C$.
The second is by the commutativity of tensor products and coproducts.

\begin{align*}
\coprod I_{\lambda} & \cong \coprod C\otimes_R \Hom_R(C,I_{\lambda}) \\
& \cong C\otimes_R \coprod (\Hom_R(C,I_{\lambda})) \\
& \cong C\otimes_R \Hom_R(C,J) \\
& \cong J.
\end{align*}
In particular, $\coprod I_{\lambda}$ is injective.
\end{proof}

Finally, we provide an example in which this technique does not seem
to provide a straightforward way in which to extend a classical
result.  Consider the following:

\begin{Quest}
Let $(R,\m,k)$ be a local, Cohen-Macaulay ring admitting a dualizing
module $D$.  Let $C$ be a semidualizing $R$-module.  If there exists
an $R$-module $M$ of finite depth with finite $\PP_C$-projective
dimension and finite $\I_C$-injective dimension, must $R$ be
Gorenstein?
\end{Quest}

Note that, as we try to apply the aforementioned techniques, we see
that $\Hom_R(C,M)$ has finite projective dimension, while
$C\otimes_R M$ has finite injective dimension.  To apply the
classical result, we need a single module that has both finite
projective dimension and finite injective dimension.

\section*{Acknowledgments}The authors would like to thank Lars Winther
Christensen, Sean Sather-Wagstaff and the referee for helpful
comments and suggestions.



\end{document}